\input amstex
\documentstyle{amsppt}
\magnification 1200
\vcorrection{-9mm}
\input epsf
\topmatter
\title     Separating semigroup of genus 4 curves
\endtitle
\author    S.~Yu.~Orevkov
\endauthor

\def\N{\Bbb N}
\def\Z{\Bbb Z}
\def\R{\Bbb R}

\def\P{\Bbb P}

\def\RP{\Bbb{RP}}

\def\Sep{\text{\rm Sep}}
\def\conj{\text{\rm conj}}
\def\supp{\text{\rm supp}}
\def\Aut{\text{\rm Aut}}
\def\ga{c}

\abstract
A rational function on a real algebraic curve $C$ is called separating if
it takes real values only at real points. Such a function defines
a covering $\Bbb R C\to\Bbb{RP}^1$. Let $\ga_1,\dots,\ga_r$ be connected
components of $\Bbb R C$.
M.~Kummer and K.~Shaw defined the separating semigroup
of $C$ as the set of all sequences $(d_1(f),\dots,d_r(f))$
where $f$ is a separating function and $d_i(f)$ is the
degree of the restriction of $f$ to $\ga_i$.

In the present paper we describe the separating semigroups of all genus 4 curves.
For the proofs we consider the canonical embedding of $C$ into a
quadric $X$ in $\P^3$ and apply Abel's theorem to 1-forms on $C$ obtained as
Poincar\'e residues of certain meromorphic 2-forms. 
\endabstract

\address
Steklov Mathematical Institute, Gubkina 8, Moscow, Russia
\endaddress

\address
IMT, l'universit\'e Paul Sabatier, 118 route de Narbonne, Toulouse, France
\endaddress

\email
orevkov\@math.ups-tlse.fr
\endemail

\thanks
A part of the work on the paper was done during author's visit to Geneva supported by the Swiss National Science Foundation project 200400
\endthanks

\endtopmatter

\def\sectLem {2}
\def\sectProof {3}
\def\sectE     {3.1}
\def\sectC     {3.2}
\def\sectH     {3.3}
\def\sectHE {4}

\def\thGenusFour {1}

\def\lemGAFA     {2.1}
\def\lemDeform   {2.2}
\def\lemExist    {2.3}
\def\lemOval     {2.4}

\def\claimCa  {3.1}
\def\claimCb  {3.2}
\def\claimCc  {3.3}
\def\claimCd  {3.4}
\def\claimCe  {3.5}
\def\claimCf  {3.6}

\def\lemH     {3.7}
\def\lemHL    {3.8}

\def\claimHa  {3.9}
\def\claimHb  {3.10}
\def\claimHc  {3.11}
\def\claimHd  {3.12}

\def\thHE  {2}

\def\eqLem   {1}
\def\eqLemC  {2}
\def\eqAB    {3}
\def\eqEC    {4}
\def\eqHE    {5}

\def\figHyp {1}
\def\figLem {2}
\def\figHE  {3}

\def\refDZ   {1}
\def\refKS   {2}
\def\refSep  {3}
\def\refGAFA {4}
\def\refZ    {5}

\document

\head 1. Introduction
\endhead

By a {\it real algebraic curve} we mean a complex algebraic curve $C$
endowed with an antiholomorphic involution $\conj:C\to C$
(the complex conjugation involution). In this case
we denote the {\it real locus} $\{p\in C\mid \conj(p)=p\}$ by $\R C$.
A real curve is called {\it dividing} or {\it separating} if $\R C$
divides $C$ into two halves exchanged by the complex conjugation.
All curves considered here are smooth and irreducible.

A necessary and sufficient condition for $C$ to be separating is the existence
of a {\it separating morphism} $f:C\to\P^1$, that is a morphism such that
$f^{-1}(\RP^1)=\R C$.
The restriction of a separating morphism to $\R C$ is a covering over $\RP^1$.
If we fix a numbering of the connected components $\ga_1,\dots,\ga_r$ of $\R C$,
we may consider the sequence $d(f)=(d_1,\dots,d_r)$ where $d_i$ is
the covering degree of $f$ restricted to $\ga_i$.
Kummer and Shaw [\refKS] defined the {\it separating semigroup} of $C$ as
$$
    \Sep(C) = \{d(f)\mid\,f:C\to \P^1 \text{\ is a separating morphism}\}.
$$
It is easy to check that this is indeed a semigroup (see [\refKS, Prop.~2.1]).
We denote:
$$
     \N=\{n\in\Z\mid n\ge 1\}, \qquad  \N_0=\{n\in\Z\mid n\ge 0\}.
$$

It is shown in [\refKS] that $\Sep(C)=\N^{g+1}$ if $C$ is an {\it $M$-curve} of genus $g$
(i.e.~$\R C$ has $g+1$ connected components) and $\Sep(C)$ is $\N^2$ (resp. $\N^3$ or $2+\N_0$) if
$C$ is a separating curve of genus 1 (resp. 2).
The separating semigroups of hyperelliptic curves of any genus and of curves of genus 3
are computed in [\refSep]. 
A much simpler proof for curves of genus 3 is given in [\refGAFA, Remark 3.3], and
in \S\sectHE\ we also give a proof for hyperelliptic curve, which is essentially the same as in [\refSep]
but exposed from the point of view proposed in [\refGAFA].

In the present paper we compute the separating semigroup of all genus 4 curves
(this result was announced in [\refGAFA, Remark 3.6]).
Let $C$ be a separating curve of genus 4. If $C$ is an $M$-curve, then $\Sep(C)=\N^5$ (see above).
If $C$ is hyperelliptic but not an $M$-curve, then $\Sep(C)=\{2\}\cup(4+\N_0)$ (see [\refSep] and \S\sectHE).

Assume now that $C$ is not hyperelliptic.
It is well-known that the image of $C$ under the canonical embedding is a degree 6 curve
on an irreducible quadric surface $X$ in $\P^3$. Since $C$ is real and separating, the
real structure on $\P^3$ can be chosen so that $C$ is a real curve on an irreducible real
surface $X$ such that $\dim\R X=2$, thus $\R X$ is an ellipsoid, a hyperboloid,
or a quadratic cone.

When $X$ is an ellipsoid or a hyperboloid, all {\it rigid isotopy classes} of
smooth real sextic curves $C$ of genus 4 on $X$ (i.e. the connected components of
the space of such curves) are described in [\refDZ]. The same arguments can be easily
adapted to the case when $X$ is a quadratic cone (see also the footnote in
[\refZ, p.~14]).
Representatives of all the rigid isotopy classes of separating non-maximal curves up to
automorphisms of $X$ are depicted in Figure~\figHyp. Four of them are realizable
as a small perturbation of  three plane sections. The other two are
perturbations of the union of a plane section and a section by a thin cylinder
around a line. The arrows in Figure~\figHyp\ represent {\it complex orientations}, i.e., the
boundary orientations induced from one of the halves of $C\setminus\R C$.

\midinsert
\centerline{\epsfxsize=100mm \epsfbox{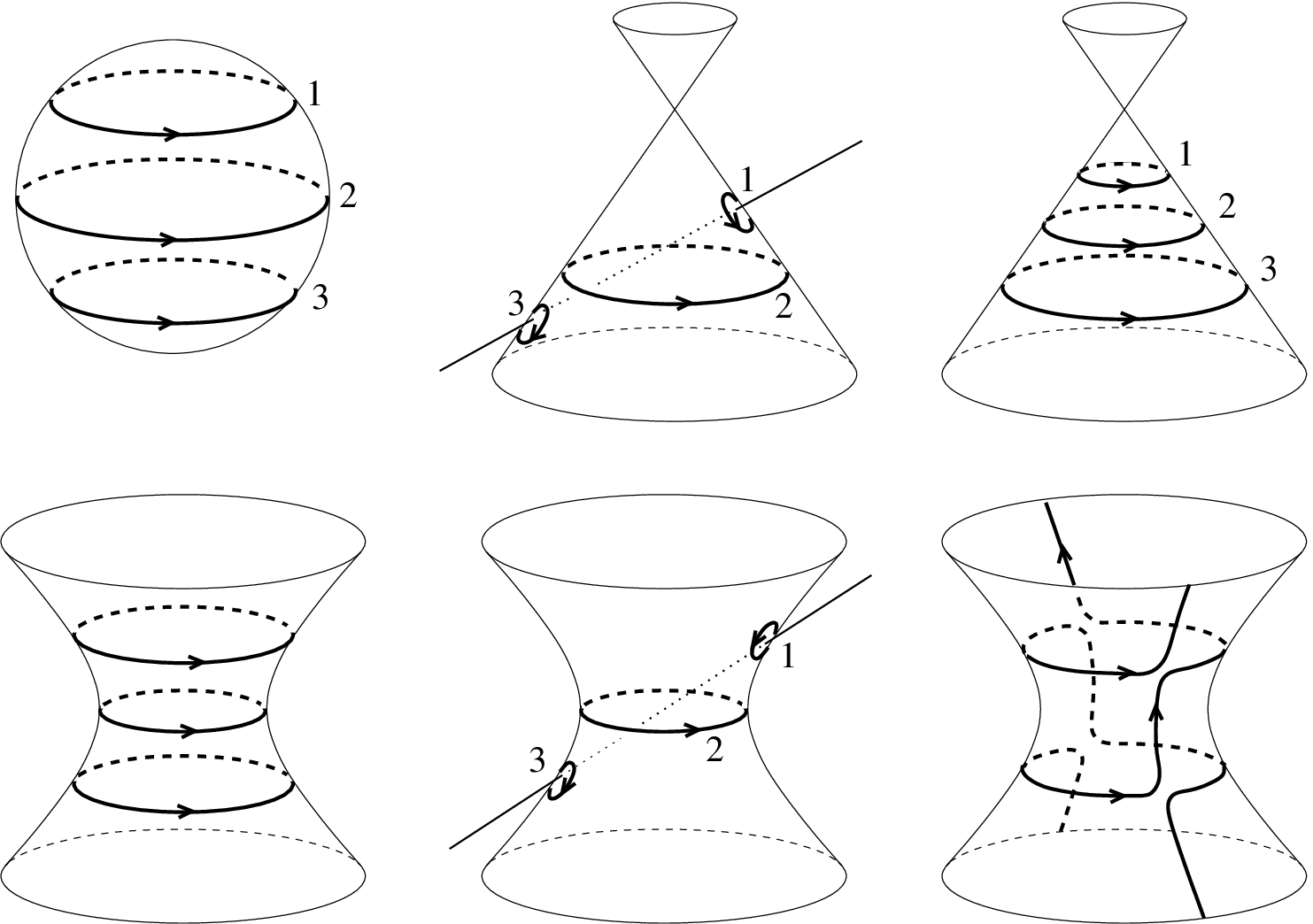}}
\botcaption{Figure~\figHyp}
    Rigid isotopy types of separating non-$M$-sextics on quadrics.
\endcaption
\endinsert

 We see in Figure~\figHyp\ that the rigid isotopy type (up to $\Aut(X)$) of a separating curve
 is determined by the topology of the pair $(\R X,\R C)$. Moreover, it is determined
 by the number of connected components of $\R C$ and the number
 of those of them which bound a smooth disk in $\R X$ (such components are called
 {\it ovals} of $C$).
 Though it is not evident a priori, it turns out that
 $\Sep(C)$ also depends on $X$ and on these two numbers only.
The main result of the paper is the following.

\proclaim{ Theorem \thGenusFour }
Let $C$ be a separating curve of genus $4$ on a real quadric $X$.
We number the components according to Figure~\figHyp\ (arbitrarily if $C$ is an $M$-curve).
Then $\Sep(C)$ is as in Table~1.
\endproclaim

\midinsert
\centerline{
\vbox{\offinterlineskip
\def\h {height2pt&\omit&&\omit&&\omit&&\omit&\cr}
    
\def\1{\bar1}\def\2{\bar2}\def\3{\bar3}\def\4{\bar4}\def\*{${}^*$}

\hrule
\halign{&\vrule#&\strut\;\;\hfil#\hfil\,\;\cr
\h
& $X$       &&$b_0(\R C)$&&number
                          of ovals&& $\Sep(C)$                 &\cr\h
 \noalign{\hrule}\h\h
& ellipsoid &&    3      &&   3   && $(1,2,1)+\N_0^3$          &\cr\h
&           &&    5      &&   5   && $\N^5$                    &\cr\h
 \noalign{\hrule}\h\h
& quadratic
  cone      &&    3      &&   0   && $\{(1,1,1)\}\cup
                                     \big((1,2,1)+\N_0^3\big)$ &\cr\h
&           &&    3      &&   2   && $(1,2,1)+\N_0^3$          &\cr\h
&           &&    5      &&   4   && $\N^5$                    &\cr\h
 \noalign{\hrule}\h\h
&hyperboloid&&    1      &&   0   && $3 + \N_0$                &\cr\h
&           &&    3      &&   0   && $\N^3$                    &\cr\h
&           &&    3      &&   2   && $(1,2,1)+\N_0^3$          &\cr\h
&           &&    5      &&   4   && $\N^5$                    &\cr\h
\noalign{\hrule}
}}
}
\botcaption{Table 1} Separating semigroups of genus 4 curves on quadrics.
\endcaption
\endinsert

Sections \sectLem\ and \sectProof\ are
devoted to the proof of Theorem~\thGenusFour.
It is based on the techniques proposed in [\refGAFA, \S3].
In Section~\sectHE\ we give a proof of [\refSep, Theorem~2] similar to the proof
of Theorem \thGenusFour\ of the present paper.

Till the end of Section~\sectProof, $X$ and $C$ are as in Theorem \thGenusFour\ and $C$ is
not an $M$-curve (Theorem \thGenusFour\ for $M$-curves is proven in
[\refKS, Thm.~1.7]). We denote the
number of connected components of $\R C$ by $r$ and the number of ovals of $C$
by $l$. As in the definition of $\Sep(C)$, the connected components of $\R C$
are denoted by $\ga_1,\dots,\ga_r$ (numbered according to Figure~\figHyp).


\head 2. Main lemmas
\endhead

Let $D=D_0+D_1$ be a real plane section of $X$ where $D_0$ is the component of $D$ of
even multiplicity (in our case $D_0$ is non-empty only when $X$ is the quadratic cone
and $D=D_0$ is a double generator).
The divisor $D-C$ belongs to the canonical class of $X$, thus it is the divisor
of some meromorphic 2-form $\Omega_D$ on $X$. It defines a ``chess-board'' orientation
on $\R X\setminus(\R C\cup \R D_1)$, i.e., an orientation which changes when crossing
$\R C\cup \R D_1$.
We define the {\it $D$-orientation} of $\R C\setminus\R D$ as the boundary orientation
induced by the ``chess-board'' orientation of $\R X\setminus(\R C\cup \R D_1)$.
Let $\omega_D$ be the Poincar\'e residue of $\Omega_D$.
Then the $D$-orientation coincides with the orientation defined by $\omega_D$
in the sense that $\omega_D(v)>0$ for $v\in T(\R C)$ if and only if the $D$-orientation is positive on $v$.
Similarly to the complex orientations, the $D$-orientation is defined up to simultaneous
reversing on all components of $\R C$.

\proclaim{ Lemma \lemGAFA } {\rm(See [\refGAFA, Thm.~3.2].)}
Let $f:C\to\P^1$ be a separating morphism and let $P=f^{-1}(x)$, $x\in\RP^1$.
If $P\not\subset D$, then the $D$-orientation cannot coincide with the
complex orientation at all points of $P\setminus D$.
\qed\endproclaim

\proclaim{ Lemma \lemDeform }
Let $p_1,\dots,p_n$ and $q_1,\dots,q_n$ be pairwise distinct points
on $C$. Suppose that the divisors
$P=p_1+\dots+p_n$ and $Q=q_1+\dots+q_n$ are linearly equivalent.
Then there exist smooth paths $p_j:[0,t_0]\to C$, $p_j(0)=p_j$, $t_0>0$, such that
$P_t := \sum p_j(t)\in|P|$ for each $t\in[0,t_0]$, and the derivative $p'_j(0)$
is nonzero for each $j=1,\dots,n$.
\endproclaim

\demo{ Proof } Let $f$ be a meromorphic function on $C$ with simple zeros at $P$ and simple poles at $Q$.
Let $P_t$ be the divisor of zeros of $f-t$. Then the result follows from the implicit function theorem.
\qed\enddemo

\proclaim{ Lemma \lemExist }
Let $D$ be an irreducible real plane section of $X$.
We fix a complex orientation on $C$ (the arrows in Figure~\figLem).

If $r=3$ and $\R D$ crosses the components of $\R C$ as shown in
Figure~\figLem\ on the left or in the middle, then
$$
    \Big((1,3,1)+\N_0^3\Big)\cup\Big((1,2,2)+\N_0^3\Big)\subset\Sep(C).   \eqno(\eqLem)
$$

If $r=1$ and $\R D$ crosses the components of $\R C$ as
in Figure~\figLem\ on the right, then $5+\N_0\subset\Sep(C)$.
\endproclaim

\midinsert
\centerline{\epsfxsize=125mm \epsfbox{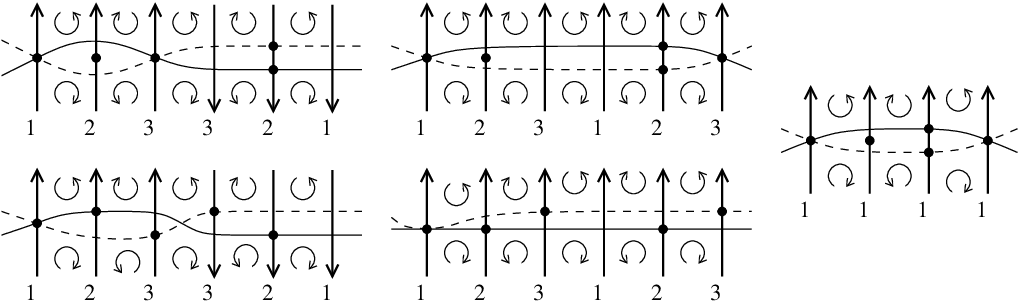}}
\botcaption{Figure~\figLem} (See Lemma~\lemExist) $\R D$ is shown by the horizontal solid line.
\endcaption
\endinsert

\demo{ Proof }
Let $D'$ be a real plane section of $X$, close to $D$ and intersecting $D\cup C$ as
shown in Figure~\figLem\ by the dashed line.
Consider the divisor $P=p_1+\dots+p_5$ on $C$, where the five points $p_1,\dots,p_5$
are placed as shown in the respective cases in Figure~\figLem.
Choosing $D'$ sufficiently close to $D$, we may assume that no line lying on $X$
passes through two points of $P$.

By Riemann-Roch Theorem the dimension of the linear system $|P|$ is at least $1$.
Let us show that $|P|$ does not have base points, which means that $\dim|P-p_j|=0$ for
each $j=1,\dots,5$. By Riemann-Roch Theorem it is enough to show that the divisor
$P-p_j$ is non-special, i.e. $|K_C-(P-p_j)|=\varnothing$.
The linear system $|K_C|$ is cut on $C$ by plane sections
(recall that $C$ is canonically embedded in $\P^3$), hence we have to show that
no plane section passes through all points of $P-p_j$.
Suppose that a plane section $D''$ passes through all points of $P-p_j$.
We consider three cases.

\smallskip
{\it Case 1.} $p_j\not\in D\cap D'$. Then $D''$ has three common points with $D$ or $D'$,
hence it coincides with $D$ or $D'$ which is impossible because each of $D$, $D'$
contains only three points of $P$.

\smallskip
{\it Case 2.} $p_j\in D\cap D'$ and $P\not\subset D\cup D'$. Then the single
point of $P\setminus(D\cup D')$ can
be moved out of the union of all the plane sections passing through
three-point subsets of $P-p_j$.

\smallskip
{\it Case 3.} $p_j\in D\cap D'$ and $P\subset D\cup D'$ (see the two bottom pictures
in Figure~\figLem). Each of $D$, $D'$ passes through two points of $P-p_j$, hence
$D''$ crosses $D\cup D'$ only at points of $P$ because $D\cdot D''=D'\cdot D''=2$.
Choosing $D'$ sufficiently close to $D$, we may ensure that both $D'$ and $D''$
are $\Cal C^1$-close to $D$, in particular, $D''$ goes monotonically from the left
to the right (according to Figure~\figLem). Easy to see that this is impossible
under condition that $D''$ passes through all points of $P-p_j$ and crosses
$D\cup D'$ only at the points of $P$.

\smallskip
We have proven that the linear system $|P|$ does not
have base points. Then by Lemma~\lemDeform\ there exists a smooth deformation
$P_t = p_1(t)+\dots+p_5(t)\in|P|$, $t\in[0,t_0]$ such that $p_j(0)=p_j$ and
$p'_j(0)\ne 0$ for each $j=1,\dots,5$.

In the rest of the proof we use the arguments as in the proof of [\refGAFA, Thm.~3.2].
We set $\omega=\omega_D$ and consider the $D$-orientation on $\R C\setminus\R D$
(see the beginning of this section). It is the boundary orientation induced
by the ``chess-board'' orientation, which is shown in Figure~\figLem\ by circular arrows.
Let $v_j=p'_j(0)\in T_{p_j}(\R C)$, $j=1,\dots,5$.
By Abel's theorem we have $\omega(v_1)+\dots+\omega(v_5)=0$
and $\omega(v_j)=0$ when $p_j\in D$. There are only two points $p_{j_1}$, $p_{j_2}$
of $P$ not belonging to $D$, hence $\omega(v_{j_1})= -\omega(v_{j_2})$.
We see in Figure~\figLem\ that the $D$-orientation and the complex orientation
coincide at one of the points $p_{j_1}$, $p_{j_2}$ and they are opposite at
the other point. Hence the complex orientations are of the same sign
on $v_{j_1}$ and $v_{j_2}$.

By applying the same arguments with $D'$ instead of $D$
we conclude that in all cases with $r=3$,
the complex orientation on the vectors $v_j$ have the same sign whenever the
corresponding points $p_j$ belong to the same connected component of $\R C$.
In the case $r=1$ (the rightmost case in Figure~\figLem), these arguments
give only that the complex orientation has the same sign on the vectors
$v_2,v_3,v_4$ (here we number the points $p_j$ from the left to the right).
However, by applying the same arguments
to the plane section through $p_2,p_4,p_5$ (resp. $p_1,p_2,p_3$), we
obtain that the complex orientation has the same sign on $v_1$ and $v_3$
(resp. on $v_4$ and $v_5$).

Therefore, for $0<t\ll t_0$, the divisors $P$ and $P_t$ are {\it interlacing}
(see [\refKS, \S2.1]), i.e., each component of $\R C\setminus P$
contains exactly one point of $P_t$. Then (see [\refKS, Prop.~2.11]) the
meromorphic function on $C$ whose divisor is $P-P_t$ defines a separating
morphism $C\to\P^1$. Hence we have $\{(1,3,1),(1,2,2)\}\subset\Sep(C)$
when $r=3$, and we have $5\in\Sep(C)$ when $r=1$.
Since the divisor $P$ is non-special,
the result follows from [\refKS, Prop.~3.2].
\qed\enddemo

\proclaim{ Lemma \lemOval } Suppose that
$X$ is a quadratic cone or hyperboloid, $l=2$, and
$L$ is a real line on $X$.
Then $L$ cannot have non-empty intersection with both ovals of $C$.
\endproclaim

\demo{ Proof } Any real curve on $X$ intersects an oval at an even number of points
counting multiplicities. Thus the result follows from the fact that $L\cdot C=3$.
\qed\enddemo


\head\sectProof. Proof of Theorem \thGenusFour
\endhead

\subhead\sectE. Proof for ellipsoids \endsubhead
Let $X$ be an ellipsoid and $r=3$.
If $(d_1,d_2,d_3)\in\Sep(C)$, then
$d_2\ge 2$; see [\refGAFA, Example 3.5].

The pencil of planes through a point of $\ga_1$ and a point
of $\ga_3$ defines a separating morphism which realizes
$(1,2,1)\in\Sep(C)$.

Let us choose a point in each component of $\R X\setminus\R C$ which is
homeomorphic to a disk.
Lemma \lemExist\ applied to a plane section $D$ passing through these two points
completes the proof of Theorem~\thGenusFour\ for ellipsoids.


\subhead\sectC. Proof for quadratic cones \endsubhead
Let $X$ be a quadratic cone. Theorem~\thGenusFour\ for this case is a combination
of the following propositions.

\proclaim{ Proposition \claimCa }
If $(r,l)=(3,0)$ and $(d_1,1,d_3)\in\Sep(C)$, then $d_1=d_3=1$.
\endproclaim

\demo{Proof}
Let $f:X\to\P^1$ be a separating morphism realizing $(d_1,1,d_3)$ and $P$ a fiber of $f$.
The result follows from Lemma~\lemGAFA\ where $D$ is a double generator passing through
the point of $P\cap \ga_2$.
\qed\enddemo

\proclaim{ Proposition \claimCb }
If $(r,l)=(3,0)$, then $\{(1,1,1),(1,2,1)\}\subset\Sep(C)$.
\endproclaim

\demo{Proof}
A pencil of plane sections passing through a real line $L$ defines
a morphism $f:X\to\P^1$. If $L$ is a generator of $X$, then $f$
realizes $(1,1,1)\in\Sep(C)$.
If $L$ meets $\ga_1$ and $\ga_3$ but not $\ga_2$, then $f$ realizes $(1,2,1)\in\Sep(C)$.
\qed\enddemo

\proclaim{ Proposition \claimCc }
If $(r,l)=(3,0)$, then
$$
   \big\{(d_1,d_2,d_3)\in(1,2,1)+\N_0^3\;\big|\; d_1+d_2+d_3\ge 5\big\} \subset \Sep(C). \eqno(\eqLemC)
$$
\endproclaim

\demo{Proof}
Let $L$ be a real line close to the rotation axis of $X$ and not passing through
the apex. Let $D$ be the section of $X$ by a real plane passing through $L$.
Then (\eqLem) holds by Lemma~\lemExist, 
and we also have $(2,2,1)+\N_0^3\subset\Sep(C)$ by symmetry.
\qed\enddemo

\proclaim{ Proposition \claimCd }
If $(r,l)=(3,2)$ and $(d_1,d_2,d_3)\in\Sep(C)$, then $d_2\ge 2$.
\endproclaim

\demo{Proof}
Suppose that there exists a separating morphism $f:X\to\P^1$ realizing $(d_1,1,d_3)$.
Let $P$ be a fiber of $f$ and let $D=2L$ where $L$ is the generator passing through
the unique point of $P\cap \ga_2$. By Lemma~\lemOval\ $L$ cannot pass through all
points of $P$. Hence we obtain a contradiction with Lemma~\lemGAFA.
\qed\enddemo

\proclaim{ Proposition \claimCe }
If $(r,l)=(3,2)$, then $(1,2,1)\in\Sep(C)$.
\endproclaim

\demo{Proof}
Suppose that $C$ is as in Figure~\figHyp, i.e. it is a perturbation of a plane
section and a thin cylinder whose axis is linked with $\ga_2$. Then the pencil of
plane sections passing through a line $L$ intersecting $\ga_1$ and $\ga_3$ realizes
$(1,2,1)\in\Sep(C)$.
As we pointed out in the introduction, any other curve with $(r,l)=(3,2)$ on $X$ is
obtained from this model curve by a continuous deformation. Such a deformation can 
be followed by a simultaneous continuous deformation of the line $L$ intersecting $\ga_1$ and $\ga_3$.
During the deformation, $L$ cannot become a generator of the cone by Lemma~\lemOval.
Thus the line $L$ remains to be linked with $\ga_2$,
hence the pencil of planes through $L$ always defines the same element of $\Sep(C)$.
\qed\enddemo

\proclaim{ Proposition \claimCf }
If $(r,l)=(3,2)$, then (\eqLemC) holds.
\endproclaim

\demo{Proof}
Let $D$ be the section of $X$ by a real plane avoiding the apex and passing through the line shown
in Figure~\figHyp. Then Lemma~\lemExist\ implies (\eqLem) and, by symmetry, (\eqLemC).
As in the proof of Proposition~\claimCe, the order of crossings of $D$ with $C$ cannot
change during a continuous deformation.
\qed\enddemo


\subhead\sectH. Proof for hyperboloids \endsubhead
Let $X$ be a hyperboloid and
let $A$ and $B$ be real lines on $X$ from different rulings.
Fix some orientations on $\R A$ and $\R B$ and denote the
corresponding homology classes in $H_1(\R X)$ by $a$ and $b$ respectively.
Then there are two rigid isotopy classes of irreducible plane sections determined by their
homology classes, which are $a+b$ and $a-b$. We assume that $A$, $B$, and
the orientations of $\R X$, $\R A$, $\R B$ are chosen so that
the horizontal and vertical plane sections in Figure~\figHyp\
(oriented according to the arrows)
belong to the classes $a+b$ and $a-b$ respectively and $ab=-ba=1$.
Thus the class of $\R C$ in $H_1(\R X)$ is
$$
             \cases
                     3a+3b, &\text{if $(r,l)=(3,3)$,}\\
                       a+b, &\text{if $(r,l)=(3,1)$,}\\
                      3a+b, &\text{if $(r,l)=(1,1)$.}
             \endcases                                        \eqno(\eqAB)
$$

The following two lemmas are easy and we omit the proofs.

\proclaim{ Lemma \lemH }
Let $p,q:[0,1]\to\R X$, $t\mapsto p_t$, $t\mapsto q_t$, be two continuous paths
such that the line $p_t q_t$ is not contained in $X$ for each $t\in[0,1]$.
Let $D_0$ be an irreducible real plane section of $X$ passing through $p_0$ and $q_0$.
Then there exists a continuous family of irreducible real plane sections $\{D_t\}_{t\in[0,1]}$
such that $\R D_t$ is homologous to $\R D_0$ and passes through $p_t$ and $q_t$ for each $t\in[0,1]$.
\qed\endproclaim

\proclaim{ Lemma \lemHL }
Let $D$ be a real irreducible plane section of $X$ such that $[\R D]=a-b$.
Let $\Gamma$ be an oriented simple closed curve on $\R X$ which belongs
to the homology class $a+b$ and has two intersection points with $\R D$.
Let $L$ be a real line passing through two points $p,q\in\R D\setminus\Gamma$.
Then $\R L$ is linked with $\Gamma$ in $\RP^3$ if and only if $\R L\cap\R D$
and $\Gamma\cap\R D$ are not interlacing on $\R D$, i.e., if and only if
$p$ and $q$ belong to the same component of $\R D\setminus\Gamma$.
\qed\endproclaim

We split Theorem~\thGenusFour\ for hyperboloids into four
Propositions~\claimHa--\claimHd\ below.
It is well-known that $X$ is biregularly equivalent to $A\times B$.
Let $\pi_A:X\to A$ and $\pi_B:X\to B$ be the projections coming from
this equivalence.

\proclaim{ Proposition \claimHa}
If $(r,l)=(1,0)$, then $\Sep(C)=3+\N_0$.
\endproclaim

\demo{Proof}
$\pi_A:\R X\to\R A$ is a 3-fold covering (see (\eqAB)),
hence it realizes $3\in\Sep(C)$.

In contrary, $\pi_B:\R X\to\R B$ is not a 3-fold covering. Indeed,
otherwise $\pi_B$ would be also a separating morphism, which contradicts (\eqAB).
Hence there exists a fiber of $\pi_B$ which has only one real intersection with $C$.
Without loss of generality we may assume that it is $A$.
Let $D$ be a small real perturbation of $A\cup B$ such that
$[\R D]=a-b\in H_1(\R X)$. Then $D$ and $C$ have $4$ real intersection points.
Let $p$ and $\bar p$ be the remaining imaginary intersections.
Consider the pencil $\Cal D$ of plane sections passing through $p$ and $\bar p$.
The real loci of its members are pairwise disjoint. Hence
any real member of $\Cal D$ is irreducible because $D$ has real
intersections with any real line in $X$.
Therefore for each $E\in\Cal D$ we have $[\R E]=[\R D]=a-b\in H_1(\R X)$
and hence $E$ has $4$ real intersections with $C$ by (\eqAB) because
$$
   \R E\cdot\R C = (a-b)(3a+b) = ab-3ba = 4ab = 4.           \eqno(\eqEC)
$$
Thus $\Cal D$ realizes $4\in\Sep(C)$.

Finally, $5+\N_0\subset\Sep(C)$ by Lemma \lemExist\ applied to any element of $\Cal D$,
because all the four intersections are positive by (\eqEC), and hence the rightmost case in
Figure~\figLem\ takes place.
\qed\enddemo

\proclaim{ Proposition \claimHb }
If $(r,l)=(3,0)$, then $\Sep(C)=\N^3$.
\endproclaim

\demo{Proof}
The projection $\pi_A:X\to A$ realizes $(1,1,1)\in\Sep(C)$.

Let $D$ be a real irreducible plane section of $X$ such that $[\R D]=a-b$.
By (\eqAB) we have $\R D\cdot\R C = 3(a-b)(a+b)=6ab$, hence $\R D$ intersects $\R C$
transversally at six points. Moreover, $\R D$ crosses the components of $\R C$ in
the order shown in Figure~\figLem\ in the middle. Indeed, the order cannot change during
a rigid isotopy, hence it is enough to check this fact for the model curve in Figure~\figHyp.
Thus (\eqLem) follows from Lemma~\lemExist. Since this result is invariant under renumbering
of components of $\R C$, we conclude that $(d_1,d_2,d_3)\in\Sep(C)$ whenever $d_1+d_2+d_3\ge 5$.

Let $p_k\in\R D\cap \ga_k$ and $p_{k+1}\in\R D\cap \ga_{k+1}$ be two consecutive points of
$\R D\cap\R C$ with respect to the order along $\R D$ (the subscripts are considered mod 3).
By Lemma~\lemHL\ the line passing through $p_k$ and $p_{k+1}$ is linked with $\ga_{k+2}$.
Hence any real plane containing this line has six real intersections with $C$.
Thus the pencil of plane sections passing through $p_k$ and $p_{k+1}$ realizes the element
$(d_1,d_2,d_3)$ of $\Sep(C)$ such that $d_k=d_{k+1}=1$ and $d_{k+2}=2$.
Thus $(d_1,d_2,d_3)\in\Sep(C)$ whenever $d_1+d_2+d_3=4$.
\qed\enddemo

\proclaim{ Proposition \claimHc }
If $(r,l)=(3,2)$ and $(d_1,d_2,d_3)\in\Sep(C)$, then $d_2\ge 2$.
\endproclaim

\demo{Proof}
Suppose that there exists a separating morphism $f:X\to\P^1$ realizing $(d_1,1,d_3)$.
The projection $\pi_A:\R C\to\R A$ is not a covering. Hence there exists a fiber (we may
assume that it is $A$) which has one real intersection point with $C$. Let $A\cap C=\{p\}$.
Then $p\in\ga_2$, hence $A\cap(\ga_1\cup\ga_3)=\varnothing$.
Let $P=f^{-1}(f(p))$ and $D=A+B$. Then $D$ cannot pass through all points of $P$ by Lemma~\lemOval,
and we obtain a contradiction with Lemma~\lemGAFA.
\qed\enddemo

\proclaim{ Proposition \claimHd }
If $(r,l)=(3,2)$, then $(1,2,1)+\N_0^3\subset\Sep(C)$.
\endproclaim

\demo{Proof}
We have $(1,2,1)\in\Sep(C)$ by the same arguments as in the proof of Proposition~\claimCe.

Let us show that the inclusion (\eqLemC) holds.
The proof is also almost the same as for Proposition~\claimCf\ but we also need Lemma~\lemH.
Namely, consider a rigid isotopy $\{C_t\}_{t\in[0,1]}$ such that $C_0$ is the model curve
shown in Figure~\figHyp\ and $C_1=C$. We denote the components of $C_t$
by $\ga_{t,1}$, $\ga_{t,2}$, $\ga_{t,3}$ according to Figure~\figHyp.

Let $D_0$ be the section of $X$ by a real plane passing through the line shown in Figure~\figHyp\
and such that $[\R D_0]=[\ga_2]$ in $H_1(\R X)$.
Then Lemma~\lemExist\ implies (\eqLem) and, by symmetry, (\eqLemC) for the model curve $C_0$. 
Let us choose continuous paths $\{p_t\}$ and $\{q_t\}$ so that $p_0,q_0\in\R D_0$ and $p_t\in\ga_{t,1}$,
$q_t\in\ga_{t,3}$ for all $t$. By Lemma~\lemOval\ the line $p_t q_t$ is not contained in $X$
for all $t$. Hence by Lemma~\lemH\ there exists a continuous family of irreducible real plane sections
$\{D_t\}$ such that $p_t,q_t\in\R D_t$ and $[\R D_t]=[\ga_2]$ for each $t$.
Then the mutual arrangement of $\R D_t$ and $\R C_t$ does not change during the deformation,
hence we may apply Lemma~\lemExist\ to $C$.
\qed\enddemo


\head\sectHE. Hyperelliptic curves
\endhead

The main results of [\refSep] can be also proved using the approach of the present paper.
One of them (a description of $\Sep(C)$ when genus$(C)=3$) is already reproved in this way in [\refGAFA].
Here we give a new proof of the other.

\proclaim{ Theorem \thHE } {\rm([\refSep, Thm.~1].)}
Let $C$ be a non-maximal hyperelliptic curve of genus $g\ge 1$. Set $m=\lfloor(g+1)/2\rfloor$. Then
$$
    \Sep(C) = \cases
               \big((1,1)\N\big) \cup \big((m,m)+\N_0^2\big), &\text{if $g$ is odd,}\\
               (2\N) \cup (g+\N_0),                    &\text{if $g$ is even.}
              \endcases
                                                                                 \eqno(\eqHE)
$$
\endproclaim

\demo{ Proof }
The curve $C$ is defined by the equation $y^2=F(x)$ where $F$ is a real polynomial
of degree $2g+2$ positive on $\R$.
It can be embedded to a real Hirzebruch surface $X$ of degree $g+1$
(the fiberwise compactification of the line bundle $\Cal O(g+1)$) so that
the hyperelliptic projection $\pi:C\to\P^1$ is the restriction of the fibration $X\to\P^1$.
The restriction $\pi|_{\R C}$ is a two-fold covering over $\RP^1$.
It is trivial if $g$ is odd and non-trivial if $g$ is even. We fix an affine chart $U\subset X$
with coordinates $(x,y)$ such that $C$ and $\pi$ are given by $y^2=F(x)$ and $(x,y)\mapsto x$ respectively.
Then $\R C\cup U$ has two connected components $c_1$ and $c_2$. Each of them cuts transversally all fibers of $\pi$.
We have $K_X+C\sim (g-1)F$ where $F$ is a fiber of $\pi$.
Denote the semigroup in the right hand side of (\eqHE) by $S$.

Let us show that $\Sep(C)\subset S$. Let $f:C\to\P^1$ be a separating morphism and $P$
is its fiber contained in $U$. Let $P_i=P\cap c_i$ and $n_i=\# P_i$, $i=1,2$.
Without loss of generality we may assume that $n_1\le n_2$.
Recall that $m=\lfloor(g+1)/2\rfloor$. Suppose that $n_1<m$. Let $D_0$ be the union of $m-1$ fibers of $\pi$
such that $P_1\subset D_0$ and $D_0\setminus\pi^{-1}(P_1)$ is disjoint from $P$.
If $g$ is odd, we set $D=2D_0$; if $g$ is even, we set $D=2D_0+F$,
where $F$ is the fiber at infinity, i.e., the fiber of $\pi$ not contained in $U$. In both cases
we have $D-C\sim K_X$ and we may consider the corresponding meromorphic 2-form $\Omega_D$ on $X$ and its
Poincar\'e residue $\omega_D$ on $C$, which is a holomorphic 1-form; $\omega_D$ defines a $D$-orientation
on $\R C\setminus\supp D$ (cf.~\S\sectLem). Then $\omega_D$ vanishes
at the points of $P_1$ and the $D$-orientation coincides with the complex orientation on
$c_2\setminus\supp D$. Hence $P_2\subset D_0$ by [\refGAFA, Thm.~3.2], and hence
$P_2\subset\pi^{-1}(P_1)$. By symmetry we also have $P_1\subset\pi^{-1}(P_2)$ and hence $n_1=n_2$,
which implies $\Sep(C)\subset S$.

\midinsert
\centerline{\epsfxsize=45mm \epsfbox{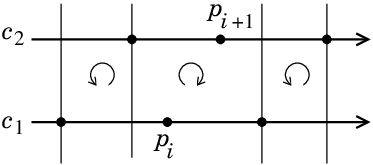}}
\botcaption{Figure~\figHE}
    To the proof of Theorem~\thHE.
\endcaption
\endinsert

Now let us prove the inverse inclusion $S\subset\Sep(C)$. The hyperelliptic projection realizes $2\in\Sep(C)$ or $(1,1)\in\Sep(C)$
(depending on the parity of $g$), hence the semigroup generated by this element is contained in $\Sep(C)$.
Thus, by [\refKS, Prop.~3.2], it is enough to realize $g\in\Sep(C)$ or $(m,m)\in\Sep(C)$ (depending on the parity of $g$)
by a separating morphism with non-special fibers.
Let $p_i=(x_i,y_i)$, $i=1,\dots,g$, be points on $C$ such that $x_1<x_2<\dots<x_g$ and sign$(y_i)=(-1)^i$,
i.e., $p_1,p_3,\dots$ are on $c_1$ and $p_2,p_4,\dots$ are on $c_2$.
Then the divisor $P=p_1+\dots+p_g$ is non-special on $C$. Let us show
that $P$ is realizable as a fiber of a separating morphism. We proceed as in the proof of Lemma~\lemExist.
By Riemann--Roch Theorem, $\dim|P|=1$. Let $P_t=p_1(t)+\dots+p_g(t)$, $P_0=P$, be a deformation of $P$ in $|P|$.
The linear system $|P|$ does not have base points. Indeed, if $p_i(t)$ is constant, then
$P^*_t:=P_t-p_i(t)\sim P^*:=P-p_i$, hence $P^*_t + \tau(P^*)\sim P^*+\tau(P^*)$, where $\tau$ is the hyperelliptic
involution, but this contradicts the fact that $P^*+\tau(P^*)\in|K_C|$ and that each element of $|K_C|$ is invariant under $\tau$.
Let $D_i$, $i=1,\dots,g-1$, be the union of $g-2$ fibers of $\pi$ passing through
all points of $P$ except $p_i$ and $p_{i+1}$. Then $D_i\in|K_X+C|$ and we consider the $D_i$-orientation on
$C\setminus\supp(D_i)$. It coincides with the complex orientation at one of the points $p_i$, $p_{i+1}$ and it
is opposite at the other point (see Figure~\figHE). Hence (cf.~the end of the proof of Lemma~\lemExist),
the complex orientation has the same sign on the tangent vectors $p'_i(0)$ and $p'_{i+1}(0)$.
This is true for all $i=1,\dots,g-1$, hence the
divisors $P$ and $P_t$, $0<t\ll 1$, are interlacing and the result follows from [\refKS, Prop.~3.2].
\enddemo



\Refs

\ref\no\refDZ\by A.\,I.~Degtyarev, V.\,I.~Zvonilov \paper Rigid isotopy
classification of real algebraic curves of bidegree $(3,3)$ on quadrics
\jour Mat. Zametki \vol 66:6 \yr 1999 \pages 810--815 \lang Russian
\transl English transl. \jour Math. Notes \vol 66 \yr 1999 \pages 670--674
\endref

\ref\no\refKS\by M.~Kummer, K.~Shaw
\paper The separating semigroup of a real curve
\jour  Ann. Fac. Sci. de Toulouse. Math\'ematiques (6) \vol 29 \yr 2020 \pages 79--96
\endref

\ref\no\refSep\by S.\,Yu.~Orevkov \paper Separating semigroup of hyperelliptic
curves and of genus 3 curves \jour Algebra i Analiz \vol 31 \yr 2019 \issue 1
\pages 108--113 \lang Russian \transl English transl. \jour St. Petersburg
Math. J. \vol 31 \yr 2020 \pages 81--84 \endref

\ref\no\refGAFA\by S.\,Yu.~Orevkov \paper Algebraically unrealizable complex
orientations of plane real pseudoholomorphic curves \jour GAFA -- Geom. Funct.
Anal. \vol 31 \yr 2021 \pages 930--947 \endref

\ref\no\refZ\by V.\,I.~Zvonilov \paper Graphs of trigonal curves and rigid
isotopies of singular real algebraic curves of bidegree $(4,3)$ on a hyperboloid
\jour arxiv:2412.15795 \endref

\endRefs
\enddocument